\documentclass[12pt]{amsart}
\usepackage{amsmath,amssymb,amsbsy,amsfonts,latexsym,amsopn,amstext,
                                               amsxtra,euscript,amscd,bm}
\usepackage{url}
\usepackage[colorlinks,linkcolor=blue,anchorcolor=blue,citecolor=blue]{hyperref}
\usepackage{color}
\usepackage{float}

\hypersetup{breaklinks=true}

\begin{document}

\newtheorem{theorem}{Theorem}
\newtheorem{lemma}[theorem]{Lemma}
\newtheorem{claim}[theorem]{Claim}
\newtheorem{cor}[theorem]{Corollary}
\newtheorem{prop}[theorem]{Proposition}
\newtheorem{definition}{Definition}
\newtheorem{question}[theorem]{Question}
\newtheorem{remark}[theorem]{Remark}
\newcommand{\hh}{{{\mathrm h}}}

\numberwithin{equation}{section}
\numberwithin{theorem}{section}
\numberwithin{table}{section}

\def\sssum{\mathop{\sum\!\sum\!\sum}}
\def\ssum{\mathop{\sum\ldots \sum}}
\def\iint{\mathop{\int\ldots \int}}

\def\squareforqed{\hbox{\rlap{$\sqcap$}$\sqcup$}}
\def\qed{\ifmmode\squareforqed\else{\unskip\nobreak\hfil
\penalty50\hskip1em\null\nobreak\hfil\squareforqed
\parfillskip=0pt\finalhyphendemerits=0\endgraf}\fi}

\newfont{\teneufm}{eufm10}
\newfont{\seveneufm}{eufm7}
\newfont{\fiveeufm}{eufm5}
%
%
\newfam\eufmfam
     \textfont\eufmfam=\teneufm
\scriptfont\eufmfam=\seveneufm
     \scriptscriptfont\eufmfam=\fiveeufm
%
%
\def\frak#1{{\fam\eufmfam\relax#1}}

\newcommand{\bflambda}{{\boldsymbol{\lambda}}}
\newcommand{\bfmu}{{\boldsymbol{\mu}}}
\newcommand{\bfxi}{{\boldsymbol{\xi}}}
\newcommand{\bfrho}{{\boldsymbol{\rho}}}

\def\fK{Frak K}
\def\fT{Frak{T}}

\def\fA{{Frak A}}
\def\fB{{Frak B}}
\def\fC{\mathfrak{C}}

\def \balpha{\bm{\alpha}}
\def \bbeta{\bm{\beta}}
\def \bgamma{\bm{\gamma}}
\def \blambda{\bm{\lambda}}
\def \bchi{\bm{\chi}}
\def \bphi{\bm{\varphi}}
\def \bpsi{\bm{\psi}}

\def\eqref#1{(\ref{#1})}

\def\vec#1{\mathbf{#1}}


\def\cA{{\mathcal A}}
\def\cB{{\mathcal B}}
\def\cC{{\mathcal C}}
\def\cD{{\mathcal D}}
\def\cE{{\mathcal E}}
\def\cF{{\mathcal F}}
\def\cG{{\mathcal G}}
\def\cH{{\mathcal H}}
\def\cI{{\mathcal I}}
\def\cJ{{\mathcal J}}
\def\cK{{\mathcal K}}
\def\cL{{\mathcal L}}
\def\cM{{\mathcal M}}
\def\cN{{\mathcal N}}
\def\cO{{\mathcal O}}
\def\cP{{\mathcal P}}
\def\cQ{{\mathcal Q}}
\def\cR{{\mathcal R}}
\def\cS{{\mathcal S}}
\def\cT{{\mathcal T}}
\def\cU{{\mathcal U}}
\def\cV{{\mathcal V}}
\def\cW{{\mathcal W}}
\def\cX{{\mathcal X}}
\def\cY{{\mathcal Y}}
\def\cZ{{\mathcal Z}}
\newcommand{\rmod}[1]{\: \mbox{mod} \: #1}

\def\cg{{\mathcal g}}

\def\vr{\mathbf r}

\def\e{{\mathbf{\,e}}}
\def\ep{{\mathbf{\,e}}_p}
\def\em{{\mathbf{\,e}}_m}

\def\Tr{{\mathrm{Tr}}}
\def\Nm{{\mathrm{Nm}}}

 \def\SS{{\mathbf{S}}}

\def\lcm{{\mathrm{lcm}}}
\def\ord{{\mathrm{ord}}}

\def\({\left(}
\def\){\right)}
\def\fl#1{\left\lfloor#1\right\rfloor}
\def\rf#1{\left\lceil#1\right\rceil}

\def\mand{\qquad \mbox{and} \qquad}

\newcommand{\commK}[1]{\marginpar{%
\begin{color}{red}
\vskip-\baselineskip 
\raggedright\footnotesize
\itshape\hrule \smallskip K: #1\par\smallskip\hrule\end{color}}}

\newcommand{\commI}[1]{\marginpar{%
\begin{color}{magenta}
\vskip-\baselineskip 
\raggedright\footnotesize
\itshape\hrule \smallskip I: #1\par\smallskip\hrule\end{color}}}

\newcommand{\commT}[1]{\marginpar{%
\begin{color}{blue}
\vskip-\baselineskip 
\raggedright\footnotesize
\itshape\hrule \smallskip T: #1\par\smallskip\hrule\end{color}}}




\hyphenation{re-pub-lished}

\mathsurround=1pt

\def\bfdefault{b}
\overfullrule=5pt

\def \F{{\mathbb F}}
\def \K{{\mathbb K}}
\def \N{{\mathbb N}}
\def \Z{{\mathbb Z}}
\def \Q{{\mathbb Q}}
\def \R{{\mathbb R}}
\def \C{{\mathbb C}}
\def\Fp{\F_p}
\def \fp{\Fp^*}

\def\Kmn{\cK_p(m,n)}
\def\psmn{\psi_p(m,n)}
\def\SI{\cS_p(\cI)}
\def\SIJ{\cS_p(\cI,\cJ)}
\def\SAIJ{\cS_p(\cA;\cI,\cJ)}
\def\SABIJ{\cS_p(\cA,\cB;\cI,\cJ)}
\def \xbar{\overline x_p}

%

\title[Cancelations between  Kloosterman sums]{Cancellations between  Kloosterman sums  modulo a prime power
with prime arguments}

 \author[K. Liu] {Kui Liu}
\address{School of Mathematics and  Statistics, Qingdao University, No.308, Ningxia Road, Shinan, Qingdao, Shandong, 266071, P. R. China}
\email{liukui@qdu.edu.cn}

 \author[I. E. Shparlinski] {Igor E. Shparlinski}

\address{Department of Pure Mathematics, University of New South Wales,
Sydney, NSW 2052, Australia}
\email{igor.shparlinski@unsw.edu.au}

 \author[T. P. Zhang] {Tianping Zhang}

\address{School of Mathematics and Information Science, Shaanxi Normal University, Xi'an 710119 Shaanxi, P. R. China}
\email{tpzhang@snnu.edu.cn}

\begin{abstract}
We  obtain a nontrivial bound for  cancellations between the Kloosterman sums modulo a large prime power with a prime argument running over very short interval, which in turn is based on a new estimate on bilinear sums of Kloosterman sums.
These results are analogues of those obtained by various authors for Kloosterman sums modulo a prime.
However the underlying technique is different and allows us to obtain nontrivial results starting
from much shorter ranges.
 \end{abstract}


\maketitle

\section{Introduction}

\subsection{Background and motivation}
For an integer $q \ge 1$ and arbitrary integers $m$ and $n$, we define the Kloosterman sums
$$
\cK(m,n;q)=\frac{1}{q^{1/2}}\sideset{}{^*}\sum\limits_{b\bmod q}e\left(\frac{mb+n\overline{b}}{q}\right),
$$
where
$\sum^*$ means summing over reduced
residue classes,  $\overline{b}$ is given by $b\overline{b}\equiv 1\bmod q$ and for a real $x$ we denote
$$
e(x)=e^{2\pi i x}.
$$

First we recall that the Weil bound, see~\cite[Corollary~11.12]{IwKow},
 yields the estimate
\begin{equation}
\label{eq:Weil}
|\cK(m,n;q)|\le \sqrt{ \gcd(m,n,q)} q^{o(1)},
\end{equation}
which gives an optimal bound on individual Kloosterman sums. So any further progress in their
applications is possible only via studying their finer structure and possible cancellations
in various families of these sums.

Indeed, starting from the groundbreaking results of  Kuznetsov~\cite{Kuz}
towards the Linnik conjecture~\cite{Lin} on cancellations between Kloosterman sums,  this has been a very active area
of research.  The initial conjecture of Linnik~\cite{Lin} predicts cancellations when the coefficients $m$ and $n$ are fixed
but the modulus $q$ varies, see~\cite{SarTsi} for results towards a uniform version of this conjecture.

The dual question, concerning cancellations between Kloosterman sums $\cK(m,n;q)$ modulo the same integer
$q$  but with varying coefficients $m$ and $n$ when one of both runs through
a consecutive interval, see~\cite{FMRS,Nied,Shp1} and references therein.
Questions of this type usually become much harder when one imposes arithmetic conditions
on the parameters, such as square-freeness, smoothness or primality. In this direction,
in the case of prime moduli $q = p$, Blomer,   Fouvry,   Kowalski,   Michel and  Mili{\'c}evi{\'c}~\cite{BFKMM2}
have shown  the existence of nontrivial cancellations
between the Kloosterman sums
$\cK(\ell,a;p)$ when $\ell$ runs over primes up to some parameter $X \ge q^{3/4 + \varepsilon}$
 for any fixed $\varepsilon > 0$
(with the goal to have $X$ as small as possible compared to $q$).

In this paper, using a  different technique,  we obtain
analogues of these results in the case where $q = p^k$  is a power of an odd prime $p$
for some sufficiently large integer $k$, which in fact allows us to reduce the amount of averaging $X \ge q^\varepsilon$ for any fixed $\varepsilon > 0$.

As in~\cite{BFKMM2}, our approach is based on bounding certain bilinear sums
of Kloosterman sums. Namely, given two sequences of weights $\cA = \{\alpha_m\}_{m=1}^M$ and
 $\cB = \{\beta_n\}_{n=1}^N$ with
 $$
 \max_{m=1, \ldots, M} |\alpha_m| \le 1 \mand
  \max_{n=1, \ldots, N} |\beta_n| \le 1,
 $$
 we consider the bilinear sums
 $$
 S_{a,q}(\cA, \cB; M, N) = \sum_{m =1}^M \sum_{n =1}^N  \alpha_m \beta_n  \cK(mn,a;q).
 $$
 Recently a series of bounds  has been obtained, by various methods, on the  sum
 $S_{a,q}(\cA, \cB; M, N)$ when  $q=p$ is prime, see~\cite{BFKMM1,BFKMM2,KMS,Shp2,ShpZha} and references therein.
 Here, using a  different technique,  we obtain analogues of these results in the case where  $q = p^k$  is a large  power of an odd prime $p$. We believe this bound can be of independent interest and may find some other applications.

\subsection{Main results}
Before we formulate our results we need to recall that  the notations  $U \ll V$ and $U = O(V)$, are equivalent
to $|U|  \le c V)$ for some constant $c>0$.
We note that we do not number implicit  constants, which we usually
also denote by $c$ so they are allowed to change their values in different statements.

We write $\ll_\rho$ and $O_\rho$ to indicate that
this constant  may depend on the parameter $\rho$.

\begin{theorem}\label{thm:Bilinear sum}
Let $q = p^k$, $k\in \Z$ be a power of an odd prime $p$. Then for any fixed constant $0<\lambda<1$ and $q^{\lambda}<M,N<q$, there exist
a constant $H_1(\lambda)>0$ depending only on $\lambda$
and an absolute constant $c>0$ such that for any $k>H_1(\lambda)$ and $\tau(\lambda)=c\lambda^3$
we have
$$\
 S_{a,q}(\cA, \cB; M, N) \ll_\lambda MNq^{-\tau(\lambda)},
$$
uniformly  over odd primes $p$ and integers $a$ with $\gcd(a,p)=1$, where the implied constant depends only on $\lambda$.
\end{theorem}

We note that for some applications one also need to restrict the summation to
a hyperbolic domain.

\begin{cor}\label{cor:Bilinear sum Hyperb}
Let   $q = p^k$, $k\in \Z$ be a power of an odd prime $p$. Then for any fixed constant $0<\lambda<1$ and $q^{\lambda}<U,V < X$ with $X<q$, there exist
a constant $H_2(\lambda)>0$ depending only on $\lambda$
and an absolute constant $c>0$ such that for any $k>H_2(\lambda)$ and $\vartheta(\lambda)=c\lambda^3$
we have
$$
 \sum_{\substack{m \ge U, \ n \ge V\\mn \le X}}   \alpha_m \beta_n  \cK(mn,a;q)
\ll_\lambda  X q^{-\vartheta(\lambda)},
$$
uniformly  over odd primes $p$ and integers $a$ with $\gcd(a,p)=1$, 
where the implied constant depends only on $\lambda$.
\end{cor}

We use Corollary~\ref{cor:Bilinear sum Hyperb}
in a combination with the main result of~\cite{LSZ} to derive:

\begin{theorem}\label{thm:Sum over primes}
Let  $q = p^k$, $k\in \Z$ be a power of an odd prime $p$ and $q^{\varepsilon}<X\leq q$. Then there exist
a constant $H_3(\varepsilon)>0$ depending only on $\lambda$
and an absolute constant $c>0$ such that for any $k>H_3(\varepsilon)$ and $\delta(\varepsilon)=c\varepsilon^3$
we have
$$
\sum\limits_{\substack{\ell \leq X\\\ell~\text{prime}}}\cK(\ell,a;q)\ll_\varepsilon Xq^{-\delta(\varepsilon)}
$$
uniformly  over odd primes $p$ and integers $a$ with $\gcd(a,p)=1$, where the implied constant depends only on $\varepsilon$.
\end{theorem}

\section{Preliminaries}

\subsection{Some bounds on Kloosterman sums}

Let $\Re\, z$ denote the real part of a complex number $z$.

First we record the following explicit formula, see, for example,~\cite[Equation~(12.39)]{IwKow}.

\begin{lemma}\label{lem:Expression for Kloosterman sums with prime power moduli}
Let $q=p^k$ with $p$ being an odd prime and $k\geq 2$, $k\in\Z$. Then for $\gcd(a,q)=1$ we have
$$
\cK(n,a;q)=\left\{
\begin{array}{ll}
2\left(\frac{r}{p}\right)^k \Re\, \vartheta_q e\left(\frac{2r}{q}\right),\quad  &{\text{if}}\quad \left(\frac{na}{p}\right)=1,\\
0,\quad &{\text{if}}\quad \left(\frac{na}{p}\right)=-1,
\end{array}
\right.
$$
where $r$ is a solution of $r^2\equiv na \bmod q$, $\left(\frac{r}{p}\right)$ is the Legendre symbol, $\vartheta_q$ equals $1$ if $q\equiv1 \bmod 4$ and $i$ if $q\equiv3 \bmod 4$.
\end{lemma}

Easy calculations show the following well-know fact (see also~\cite{LSZ}):

\begin{lemma}\label{lem:Special Kloosterman sums equals to 0}
Let $p$ be an odd prime and $k\geq 2$ be a positive integer. If $\gcd(a,p)=1$ and $p \mid n$, then the Kloosterman sums $\cK(n,a;p^k)=0$.
\end{lemma}

One of our principal technical tools is the following estimate from~\cite{LSZ} on the cancellations
amongst Kloosterman sums,

\begin{lemma}\label{lem:type1 estimate}
Let $q = p^k$, $k\in \Z$ be a power of an odd prime $p$. Then for any fixed constant $0<\lambda<1$ and $q^{\lambda}<N<q$, there exist a constant $K_0(\lambda)$ depending only on $\lambda$ and an absolute constant $c>0$ such that for any $k>K_0(\lambda)$ and $\tau(\lambda)=c\lambda^3$ we have
$$
\sum\limits_{1\leq n\leq N}\cK(n,a;q)\ll_{\lambda}Nq^{-\tau(\lambda)}
$$
uniformly  over odd primes $p$ and integers $a$ with $\gcd(a,p)=1$, where the implied constant depends only on $\lambda$.
\end{lemma}

\subsection{Short exponential sums with special polynomials}

The following bound  is contained in the proof of~\cite[Theorem~1.2]{LSZ}
(and in fact is in the core of this proof).

\begin{lemma}\label{lem:bound of exponential sum}
Suppose $\gcd(h,q)=1$, $q=p^k$ with $p$ an odd prime and $k\in\Z$.  For $q^{\eta}<N<q$ with $0<\eta<1$ being a fixed constant, there exists a constant $k_0(\eta)>0$, such that for the   polynomial
by
\begin{equation}
\label{eq:f(t)}
f(X)=\sum\limits_{j=0}^{\fl{k/s}}g(j) \gamma^j
p^{js}X^j,
\end{equation}
where $\gamma$ is an arbitrary integer with $\gcd(\gamma,p)=1$,
$$
s=\fl{\frac{\eta \log N }{3000 \log p}},
$$
$g(0)=1$ and $g(j)$ with $1\leq j\leq \fl{k/s}$ are integers given by
$$
g(j)\equiv \frac{1/2(1/2-1)\cdots(1/2-j+1)}{j!}\bmod p^k,\qquad  0\leq g(j)<p^k,
$$
and any $k>k_0(\eta)$, we have
$$
\sum\limits_{n\leq N}e\left(\frac{hf(n)}{q}\right)\ll_\eta Nq^{-\rho(\eta)},
$$
where $\rho(\eta)=c\eta^3$ with some absolute constant $c>0$ and the implied constant depends only on $\eta$.
\end{lemma}

\subsection{Vaughan identity}

As usual, we use $\Lambda(n)$ to denote the  von Mangoldt function
$$\Lambda(n)=
\begin{cases}
\log \ell &\qquad\text{if $n$ is a power of a prime $\ell$,} \\
0&\qquad\text{otherwise.}
\end{cases}$$

We use the following result of Vaughan~\cite{Va} in the form given by Davenport~\cite[Chapter~24]{Dav}.

\begin{lemma}\label{lem:Vaughan's identity}
For any complex-valued function $f(n)$ and any real numbers $1 < U,V\leq X$, we have
$$
\sum\limits_{n\leq X}\Lambda(n)f(n)\ll\Sigma_1+\Sigma_2 \log(UV)+ \Sigma_3 \log X+\left|\Sigma_4\right|,
$$
where
\begin{align*}
\Sigma_1 &=\left|\sum\limits_{n\leq U}\Lambda(n)f(n)\right|,\\
\Sigma_2 & = \sum\limits_{v\leq UV}\left|\sum\limits_{s\leq N/v}f(sv)\right|,\\
\Sigma_3 &= \sum\limits_{v\leq V}\max\limits_{w\geq 1}\left|\sum\limits_{w\leq s\leq N/v}f(sv)\right|
\end{align*}
and
$$
\Sigma_4 =\sum\limits_{\substack{uv\leq X\\u>U,~v>V}}\Lambda(u)\(\sum\limits_{d\mid v,~d\leq V}\mu(d)\)f(uv).
$$
\end{lemma}

\section{Proof of main results}

\subsection{Proof of Theorem~\ref{thm:Bilinear sum}}

Note that due to the symmetry of the sum and the claimed bound with respect to $M$ and $N$, we can
always assume that $M \ge N$.

Let $1\leq s\leq k$ be a parameter to be determined. Write
$$
 S_{a,q}(\cA, \cB; M, N)=\sum\limits_{u=1}^{p^s}\sum\limits_{v=1}^{p^s}\sum\limits_{\substack{m\leq M\\m \equiv u \bmod{p^s}}}\sum\limits_{\substack{n\leq N\\n \equiv v \bmod{p^s}}}
\alpha_m\beta_n\cK(mn,a;q).
$$
Then by Lemma~\ref{lem:Special Kloosterman sums equals to 0}, we have
$$
 S_{a,q}(\cA, \cB; M, N)\le \sideset{}{^*}\sum\limits_{u=1}^{p^s}\sideset{}{^*}\sum\limits_{v=1}^{p^s}\sum\limits_{\substack{n\leq N\\n \equiv v \bmod{p^s}}}\left|\sum\limits_{\substack{m\leq M\\m \equiv u \bmod{p^s}}}
\alpha_m\cK(mn,a;q)\right|.
$$
Now by Lemma~\ref{lem:Expression for Kloosterman sums with prime power moduli}, we get
\begin{align*}
 S_{a,q}&(\cA, \cB; M, N)\\
 &\le 2\sideset{}{^*}\sum\limits_{u=1}^{p^s}\sideset{}{^*}\sum\limits_{v=1}^{p^s}\sum\limits_{\substack{n\leq N\\n \equiv v \bmod{p^s}}}\left|\sum\limits_{\substack{m\leq M\\m \equiv u \bmod{p^s}\\\left(\frac{amn}{p}\right)=1}}
\alpha_m\left(\frac{r_{m,n}}{p}\right)^k\Re\, \vartheta_q e\left(\frac{2r_{m,n}}{q}\right)\right|,
\end{align*}
where $r_{m,n}$ is a solution to $r_{m,n}^2\equiv a mn  \bmod q$ (the inequality holds for any choice
of this solution).
Thus
\begin{equation}
\label{eq:S_a,q(M,N)}
 S_{a,q}(\cA, \cB; M, N)\le 2\underset{\left(\frac{auv}{p}\right)=1}{\sideset{}{^*}\sum\limits_{u=1}^{p^s}\sideset{}{^*}\sum\limits_{v=1}^{p^s}} T_{s}(u,v),
\end{equation}
where
\begin{equation}
\label{eq:Tbc}
T_{s}(u,v)=\sum\limits_{\substack{n\leq N\\n \equiv v \bmod{p^s}}}\left|\sum\limits_{\substack{m\leq M\\m \equiv u \bmod{p^s}}}
\alpha_m\left(\frac{r_{m,n}}{p}\right)^k\Re\, \vartheta_q e\left(\frac{2r_{m,n}}{q}\right)\right|.
\end{equation}
Note that $r_{m,n}^2\equiv amn \equiv auv \bmod p$, thus the Legendre symbol $\left(\frac{r_{m,n}}{p}\right)$ does not depend on $m$ and $n$.
We now use the Cauchy inequality and derive
$$
 T_{s}(u,v)^2
\ll(Np^{-s}+1)\sum\limits_{\substack{m_1,m_2\leq M\\m_1\equiv u\bmod p^s\\m_2\equiv u\bmod p^s}}
\left|\sum\limits_{\substack{n\leq N\\n \equiv v \bmod{p^s}}}e\left(\frac{2(r_{m_1,n}-r_{m_2,n})}{q}\right)\right|,
$$
with $r_{m_i,n}$ being a solution of the congruence $r^2\equiv m_{i}na \bmod q$ for $i=1, 2$.

Now we consider the inner sum over $n$. Note that for $i=1,2$, we have $\gcd(am_i ,q)=1$, then we use the following argument, which is similar to that in~\cite{Kha} and~\cite{LSZ}.

By $n\equiv v\bmod p^s$, there exists $t\in \Z$, such that $n=v+p^st$. Now we have
$$
r_{m_i,n}^2\equiv a m_i(v+p^st)\equiv a m_iv(1+\bar{v} p^st)\ \bmod q,
$$
with $v\bar{v}\equiv1\bmod q$. Note that $a m_iv\equiv auv\bmod p$, then
$$
\(\frac{am_iv}{p}\)=1.
$$
By the {\it Hensel lifting\/}, there exists an integer $\omega$, such that $0\leq \omega< p^s$ and
$$
\omega^2\equiv auv\bmod p^s,
$$
and also for each $m_i$ there exists an integer $\omega_{m_i}$ such that
$$
\omega_{m_i}^2\equiv a m_iv \bmod q
$$
and $\omega_{m_i}\equiv \omega\bmod p^s$ for $i=1,2.$
Thus
$$
r_{m_i,n}^2\equiv a m_iv \equiv \omega_{m_i}^2(1+\bar{v}p^st)\ \bmod q.
$$
We remark that $\omega_{m_i}$ is determined by $a$, $m_i$ and $v$,  and does not depend on $n$.
Consider $1+\bar{v}p^st$ in the $p$-adic field $\Q_p$. By the Taylor expansion
(see~\cite[Chapter~IV.1]{Ko}), we have
$$
(1+\bar{v}p^st)^{1/2}=1+\sum\limits_{j=1}^{\infty}\binom{1/2}{j}\bar{v}^jp^{js}t^j
$$
for $s\geq 1$ with the coefficients
$$\binom{1/2}{j}=\frac{1/2(1/2-1)\cdots(1/2-j+1)}{j!}, \qquad j =1,2,\ldots,
$$
which are  $p$-adic integers, since $p$ is an odd prime.
Then we have
$$
(1+\bar{v}p^st)^{1/2}\equiv 1+\sum\limits_{j=1}^{\fl{k/s}}g(j)\bar{v}^jp^{js}t^j\ \bmod p^k,
$$
where $g(j)$ with $1\leq j\leq \fl{k/s}$ are integers given by
\begin{equation}
\label{eq: Expression of g(j)}
g(j)\equiv\binom{1/2}{j}\bmod p^k,\qquad  0\leq g(j)<p^k.
\end{equation}
Thus we get two solutions for the quadratic congruence of $r_{m_i,n}$ in the inner sum
in~\eqref{eq:Tbc}.
$$
r_{m_i,n}\equiv\pm\omega_{m_i} f(t)\bmod q,
$$
where the polynomial  $f(X)$ is given by~\eqref{eq:f(t)} as in Lemma~\ref{lem:bound of exponential sum}
with $\gamma = \bar{v}$.

Choosing the solution $r_{m_i,n}\equiv\omega_{m_i} f(t)\bmod q$, then we have
$$
T_{s}(u,v)^2
\leq(Np^{-s}+1)\sum\limits_{\substack{m_1,m_2\leq M\\m_1\equiv u\bmod p^s\\m_2\equiv u\bmod p^s}}
\left|\sum\limits_{t\leq (N-v)/p^s}e\left(\frac{2(\omega_{m_1}-\omega_{m_2})f(t)}{q}\right)\right|.
$$
Since $v\leq p^s$, we can replace the range of $t$ by $t\leq N/p^s$ up to a small  error term and
thus we get
\begin{equation}
\label{eq:T(u,v) square 1}
\begin{aligned}
T_{s}(&u,v)^2 \ll (Np^{-s}+1)\sum\limits_{\substack{m_1,m_2\leq M\\m_1\equiv u\bmod p^s\\m_2\equiv u\bmod p^s}}
\left|\sum\limits_{t\leq N/p^s}e\left(\frac{2(\omega_{m_1}-\omega_{m_2})f(t)}{q}\right)\right|\\
&\qquad\qquad\qquad\qquad\qquad\qquad\qquad +(Np^{-s}+1)(Mp^{-s}+1)^2.
\end{aligned}
\end{equation}
We now choose
\begin{equation}
\label{eq: Choice s}
B(\lambda)=3000/\lambda>3000 \mand s=\fl{\frac{\log N}{B(\lambda)\log p}}.
\end{equation}
 Since $q^\lambda<N \le M<q=p^k$, then we have $s\leq k/B(\lambda)$.
 Hence,
 \begin{equation}
\label{eq: Useful Ineq}
N/p^s\geq q^{\lambda/2} \mand k-3s>k/2.
\end{equation}

In particular we see from~\eqref{eq: Useful Ineq} that
$$
Mp^{-s}+1 \ll Mp^{-s}  \mand Np^{-s}+1 \ll Np^{-s}.
$$
Thus~\eqref{eq:T(u,v) square 1} can be simplified as
\begin{equation}
\label{eq:T(u,v) square 2}
\begin{aligned}
T_{s}(u,v)^2 & \ll Np^{-s}\sum\limits_{\substack{m_1,m_2\leq M\\m_1\equiv u\bmod p^s\\m_2\equiv u\bmod p^s}}
\left|\sum\limits_{t\leq N/p^s}e\left(\frac{2(\omega_{m_1}-\omega_{m_2})f(t)}{q}\right)\right|\\
&\qquad\qquad\qquad\qquad\qquad\qquad\qquad \qquad\quad + M^2Np^{-3s}.
\end{aligned}
\end{equation}

Let $\ord_p(n)$ denote the $p$-adic order of integer $n$. For a fixed $m_1$, we claim that
$$
\sharp~\{m_2\leq M~:~m_2\equiv m_1\bmod\ p^s,~\ord_p(\omega_{m_1}-\omega_{m_2})\geq 3s\}\leq Mp^{-3s}+1.
$$
To see this, we note that $\ord_p(\omega_{m_1}-\omega_{m_2})\geq 3s$ implies
$$
\omega^2_{m_1}\equiv\omega^2_{m_2} \bmod\ p^{3s}.
$$
Since $\gcd(av, q)=1$, then by the definition of $\omega_{m_i}$, we have the congruence $m_2\equiv m_1 \bmod\ p^{3s}$, which yields our claim.

We  estimate
the contribution of this case to $T_{s}(u,v)^2$ trivially as
$$
Np^{-s}\sum\limits_{\substack{m_1\leq M\\m_1\equiv u\bmod p^s}}\sum\limits_{\substack{m_2\leq M\\m_2\equiv m_1\bmod p^s\\\ord_p(\omega_{m_1}-\omega_{m_2})\geq 3s}}
\sum\limits_{t\leq N/p^s}1 \ll M N^2 p^{-3s}(Mp^{-3s}+1).
$$
Hence, we only need to estimate the exponential sum
$$
\sum\limits_{t\leq N/p^s}e\left(\frac{2(\omega_{m_1}-\omega_{m_2})f(t)}{q}\right),
$$
with $\ord_p(\omega_{m_1}-\omega_{m_2})<3s.$ After cancelling the factor $p^s$, we can get a similar exponential sum of the type
$$
\sum\limits_{t\leq N/p^s}e\left(\frac{hf(t)}{p^r}\right)
$$
with $\gcd(h, q)=1$ and $r>k-3s$. To bound this sum, we recall our choice~\eqref{eq: Choice s} and apply
Lemma~\ref{lem:bound of exponential sum} and get
$$
\sum\limits_{t\leq N/p^s}e\left(\frac{2(\omega_{m_1}-\omega_{m_2})f(t)}{q}\right)\ll Np^{-s}q^{-\rho(\lambda/2)},
$$
for $k\geq k_0(\lambda/2)$, where $k_0$ and $\rho$ are the same as in Lemma~\ref{lem:bound of exponential sum}. Insert this bound to~\eqref{eq:T(u,v) square 2}, we get
\begin{align*}
T_{s}(u,v)^2
&\ll~ M^2N^2 p^{-4s} q^{-\rho(\lambda/2)}+ M N^2 p^{-3s}(Mp^{-3s}+1) + M^2Np^{-3s}\\
&= M^2N^2 p^{-4s} q^{-\rho(\lambda/2)} + M^2 N^2 p^{-6s}  + M^2Np^{-3s}.
\end{align*}
Then we have
$$
T_{s}(u,v)\ll MNp^{-2s}q^{-\rho(\lambda/2)/2}+MNp^{-3s} +  MN^{1/2}p^{-3s/2},
$$
which after the substitution in~\eqref{eq:S_a,q(M,N)} yields
$$
 S_{a,q}(\cA, \cB; M, N)\ll MNq^{-\rho(\lambda/2)/2}+MN^{1/2}p^{s/2}+MNp^{-s}.
$$
Note that the definition of $s$ implies
that
\begin{align*}
s  & \ge \frac{\log N}{B(\lambda)\log p} - 1 =  \frac{\lambda \log N}{3000 \log p}- 1
\\
& =  \frac{\lambda k \log N}{3000 \log q}- 1 \ge  \frac{\lambda k \log N}{6000 \log q} =  \frac{\log N}{2 B(\lambda)\log p},
\end{align*}
provided that $k$ is large enough (we also recall that $q> M\ge N> q^\lambda$).
Hence
$$
 N^{1/(2B(\lambda))}  \leq p^s
\leq N^{1/B(\lambda)},
$$
from which the result follows.

\subsection{Proof of  Corollary~\ref{cor:Bilinear sum Hyperb}}
Separating the summation over $m$ in dyadic ranges and slightly reducing the
constant $c$ in the definition of $\vartheta(\lambda)$,
 it is enough to show that for any $M_1, M_2$ with
 $$
M_2 \le 2 M_1 \mand U \le M_1 < M_2 \le 2X/V,
$$
we have
\begin{equation}
\label{eq: diadic sum}
 \sum_{m =M_1}^{M_2} \sum_{V \le n \le X/m}    \alpha_m \beta_n  \cK(mn,a;q)
\ll_\lambda X q^{-\vartheta(\lambda)}.
\end{equation}
Let $L = \fl{X/M_1}$. Clearly we can assume that $L \ge V$ as otherwise the sum over $n$
is void and there is nothing to prove.
In particular, we can assume that
$$
L \ge q^\lambda.
$$
By  the orthogonality of exponential functions, we write
\begin{align*}
 \sum_{m =M_1}^{M_2}& \sum_{V \le n \le X/m}    \alpha_m \beta_n  \cK(mn,a;q)\\
& =  \sum_{m =M_1}^{M_2}  \sum_{V \le n \le L}
\frac{1}{L}  \sum_{z=0}^{L-1} \sum_{w=0}^{X/m} e(z(w-n)/L) \alpha_m \beta_n  \cK(mn,a;q)\\
& =  \frac{1}{L}  \sum_{z=0}^{L-1}  \sum_{m =M_1}^{M_2}  \sum_{w=0}^{X/m} e(z w/L)
 \sum_{V \le n \le L} \alpha_m  \beta_n e(-z n/L) \cK(mn,a;q).
\end{align*}

Using that for $z = 0, \ldots, L-1$ we have
$$
\sum_{w=0}^{X/m} e(z w/L)   \ll  \frac{L}{1 + \min\{z, L-z\}},
$$
see~\cite[Bound~(8.6)]{IwKow}. We now derive
\begin{equation}
\label{eq: Prelim}
\begin{split}
 & \sum_{m =M_1}^{M_2} \sum_{V \le n \le X/m}    \alpha_m \beta_n  \cK(mn,a;q)\\
   & \qquad =      \sum_{z=0}^{L-1}   \frac{1}{1 + \min\{z, L-z\}}
   \sum_{m =M_1}^{M_2}    \sum_{V \le n \le L}
\widetilde \alpha_{z,m}  \widetilde \beta_{z,n}  \cK(mn,a;q)
\end{split}
\end{equation}
with some weights $|\widetilde \alpha_{z,m}| \le 1$ and  $| \widetilde \beta_{z,n}|\le 1$
(in particular $ \widetilde \beta_{z,n} = \beta_n e(-z n/L) $). Thus, by Theorem~\ref{thm:Bilinear sum}
for every $z = 0, \ldots, L-1$, for the  double sum over $m$ and $n$ we have the bound
$$
   \sum_{m =M_1}^{M_2}    \sum_{V \le n \le L}
\widetilde \alpha_{z,m}  \widetilde \beta_{z,n}  \cK(mn,a;q) \ll M_2 Lq^{-\tau(\lambda)} \ll Xq^{-\tau(\lambda)},
$$
which after substitution in~\eqref{eq: Prelim} implies~\eqref{eq: diadic sum} and concludes the proof.
\subsection{Proof of Theorem~\ref{thm:Sum over primes}}

Using partial summation, one can easily see that it is enough to estimate   the sum
$$
\widetilde S(X;a,q)=\sum\limits_{n\leq X}\Lambda(n)\cK(n,a;q).
$$
We now apply  Lemma~\ref{lem:Vaughan's identity} (with $f(n) = \cK(n,a;q)$) where we take  $U=V=X^{1/3}$, for which we need to estimate the sums:
\begin{align*}
\Sigma_1 &=\left|\sum\limits_{n\leq U}\Lambda(n)\cK(n,a;q)\right|,\\
\Sigma_2 & = \sum\limits_{v\leq UV}\left|\sum\limits_{s\leq X/v}\cK(sv,a;q)\right|,\\
\Sigma_3 &= \sum\limits_{v\leq V}\max\limits_{w\geq 1}\left|\sum\limits_{w\leq s\leq X/v}\cK(sv,a;q)\right|
\end{align*}
and
$$
\Sigma_4 =\sum\limits_{\substack{uv\leq X\\u>U,~v>V}}
\Lambda(u)\(\sum\limits_{d\mid v,~d\leq V}\mu(d)\)\cK(uv,a;q).
$$

Using the Weil bound~\eqref{eq:Weil}, or  its  more precise version given in
Lemma~\ref{lem:Expression for Kloosterman sums with prime power moduli}, we have
$$
\Sigma_1\ll_\varepsilon X^{1/3}q^{\varepsilon/3}<Xq^{-\varepsilon/3}.
$$

For $\Sigma_2$, note that $v\leq UV=X^{2/3}$, then $X/v\geq X^{1/3}>q^{\varepsilon/3}$. Now Lemma~\ref{lem:type1 estimate} can be applied to the sum over $s$. Hence there exist constants $K_1(\varepsilon)>0$ and $\tau_1(\varepsilon)>0$ such that
$$
\Sigma_2 \ll_\varepsilon Xq^{-\tau_1(\varepsilon)}
$$
holds uniformly for any odd prime $p$ and $k>K_1(\varepsilon)$.

For $\Sigma_3$, write the inner sum over $s$ into
\begin{equation}\label{eq:inner sum of sum_3}
\sum\limits_{w\leq s\leq X/v}\cK(sv,a;q)=\sum\limits_{s\leq X/v}\cK(sv,a;q)-\sum\limits_{s< w}\cK(sv,a;q).
\end{equation}
The contribution of the first sum on the right side can be bounded similarly as $\Sigma_2$. For the contribution of the second sum,
If $w< q^{\varepsilon/4}$, then by the Weil bound we have
$$
\sum\limits_{s\leq w}\cK(sv,a;q)\ll_\varepsilon q^{\varepsilon/3}.
$$
If $w\ge q^{\varepsilon/4}$, then Lemma \ref{lem:type1 estimate} can be applied to the above sum. It follows from the above treatment that there exist constants $K_2(\varepsilon)>0$ and $\tau_2(\varepsilon)>0$ such that
$$
\Sigma_3 \ll_\varepsilon Xq^{-\tau_2(\varepsilon)}
$$
holds uniformly for any odd prime $p$ and $k>K_2(\varepsilon)$.

For $\Sigma_4$, by Corollary~\ref{cor:Bilinear sum Hyperb} there exist constants $K_3(\varepsilon)>0$ and $\tau_3(\varepsilon)>0$ such that
$$
\Sigma_4 \ll_\varepsilon Xq^{-\tau_3(\varepsilon)}
$$
holds uniformly for any odd prime $p$ and $k>K_3(\varepsilon)$.
It is also easy to see that for the functions $\tau_\nu(\varepsilon)$ we can take
$\tau_\nu(\varepsilon)= c_\nu \varepsilon^3$ for some absolute constants $c_\nu>0$,
for every $\nu=1,2,3$.

Now the desired bound on $\widetilde S(X;a,q)$ and thus Theorem~\ref{thm:Sum over primes} follow from the above estimates.

\end{document}